\documentclass{amsart}
\usepackage[all]{xy}
\usepackage{amssymb}
\usepackage{amsmath}
\usepackage{amsmath}

\swapnumbers
\newtheorem{thm}[subsection]{Theorem}
\newtheorem{prop}[subsection]{Proposition}
\newtheorem{lemma}[subsection]{Lemma}
\newtheorem{cor}[subsection]{Corollary}

\def\Ai{A_{\infty}}

\def\ac{{}^{\scriptstyle \textrm{!`}}}

\def\d{\succ}
\def\g{\prec}
\def\l{\dashv }
\def\r{\vdash }

\def\dd{\delta}

\def\t{\otimes}

\def\cc{\gamma}

\def\I{\mathrm{I}}

\def\Sy{\mathbb{S}}

\newcommand{\Vt}[1]{V^{\otimes #1}}

\def\RR{{\mathbb{R}}}

\def\RR{{\mathbb{R}}}

\def\PP{{\mathcal{P}}}

\def\TTT{{\mathcal{T}}}

\def\Id{\mathrm{Id }}

\def\Vect{\textsf{Vect}}

\def\id{\mathrm{ id }}
\def\I{\mathrm{ I }}

\def\ac{{}^{\scriptstyle \textrm{!`}}}

\def\KK{\mathbb{K}}

\def\arbreA{\vcenter{\xymatrix@R=3pt@C=3pt{
&& \\
&*{}\ar@{-}[ur] \ar@{-}[ul] \ar@{-}[d]     &\\
&&
}}}

\newenvironment{proo}{\begin{trivlist} \item{\emph{Proof.}}}
  {\hfill $\square$ \end{trivlist}}

\begin{document}

\author[J.-L. Loday]{Jean-Louis Loday}
\address{Institut de Recherche Math\'ematique Avanc\'ee\\
    CNRS et Universit\'e de Strasbourg\\
    7 rue R. Descartes\\
    67084 Strasbourg Cedex, France}
\email{loday@math.u-strasbg.fr}

\title{On the operad of associative algebras with derivation}

\dedicatory{Dedicated to Tornike Kadeishvili in honor of his sixtieth birthday}

\begin{abstract} We study the operad of associative algebras equipped with a derivation. We show that it is determined by polynomials in several variables and substitution. Replacing polynomials by rational functions gives an operad which is isomorphic to the operad of ``moulds''.  It provides an efficient environment for doing integro-differential calculus. Interesting variations are obtained by using formal group laws. The preceding case corresponds to the additive formal group law. We unravel the notion of homotopy associative algebra with derivation in the spirit of Kadeishvili's work.
\end{abstract}

\maketitle

\section*{Introduction} \label{S:int} To any one-dimensional formal group law we show how to  associate a type of algebras, that is an operad, whose space of $n$-ary operations is the algebra of rational functions in $n$ variables. If the formal group law is polynomial, then the spaces of polynomials in $n$ variables form a suboperad. We show that, for the additive formal group law, the polynomial operad is simply the operad $AsDer$ of associative algebras equipped with a derivation, and the rational functions operad $RatFct$ is isomorphic to the operads of ``moulds''.  Since the operad $RatFct$ contains both the derivation operator and its inverse, that is the integral operator, the whole integro-differential calculus can be written in terms of operadic calculus. The operad $RatFct$ is related to the dendriform and tridendriform operads, and several others.

In order to deal with formal group laws in noncommutative variables, one has to modify the notion of operad and work with ``pre-shuffle algebras'' introduced by M.\ Ronco in \cite{Ronco07}. The idea is essentially to forget parallel composition 
in the definition of an operad by means of the partial compositions. 

In the appendix we show that the operad $AsDer$, which is quadratic, is a Koszul operad. We describe explicitly the notion of ``homotopy associative algebra with derivation'', that is $AsDer_{\infty}$-algebra.
\bigskip

\noindent {\bf Notation.} In this paper $\KK$ is a commutative unital algebra and all modules over $\KK$ are supposed to be free. Its unit is denoted by 1. Since most of the time we think of $\KK$ as being a field we will say vector space or space for a free module over $\KK$.
The tensor product  of vector spaces over $\KK$
is denoted by $\t$. The tensor product of $n$ copies of the space $V$ is
denoted by $\Vt n$. The symmetric group $\Sy_{n}$ is acting on $\Vt n$ by place-permutation.


The polynomial algebra, resp.\ the algebra of rational functions,  on the set of variables $\{x_{1}, \ldots , x_{n}\}$ is denoted by $ \KK[x_{1}, \ldots , x_{n}]$, resp.\  $\KK( x_{1}, \ldots , x_{n})$. The unit is denoted by $1_{n}$ or $1$ if there is no confusion.

\section{Prerequisites on operads} This is a very brief introduction on algebraic operads, whose purpose is essentially to set up notations. One can consult for instance \cite{MSS}, or \cite{LV}, for more details.

\subsection{Nonsymmetric operad}\label{nsop} A \emph{nonsymmetric operad}, or ns operad for short, is a graded vector space $\{\PP_{n}\}_{n\geq 0}$ equipped with a particular element of $\PP_{1}$ denoted $\id$ and called the identity operation, and a family of linear maps
$$\cc_{i_{1}, \ldots , i_{k}} : \PP_{k}\t \PP_{i_{1}}\t \cdots \t \PP_{i_{k}}\to \PP_{n},\quad n=i_{1}+\cdots + i_{k},$$
which satisfy the following associative and unital properties. On the functor $\PP: \Vect \to \Vect$ defined as $\PP(V):=\bigoplus_{n}\PP_{n}\t V^{\t n}$ the operations $\cc_{i_{1}, \ldots , i_{k}}$ induce a transformation of functors $\cc : \PP\circ \PP \to \PP$. Then, $\cc$ is supposed to be associative. The element $\id$ can be interpreted as a morphism from the identity functor $\I$ to $\PP$ that we denote by $\iota : \I \to \PP$. Then,  $\iota$ is supposed to be a unit for $\cc$.

There is an alternative definition of a ns operad which uses the so-called partial compositions. For any $i\leq m$ the \emph{partial composition}
$$\circ_{i}:\PP_{m}\t \PP_{n} \to \PP_{m+n-1}$$
evaluated on $(\mu, \nu)$ is the operation $\cc_{1, \ldots ,1, n, 1, \ldots , 1}$ evaluated on $(\mu ; \id, \ldots \id, \nu, \id, \ldots, \id)$. 
Pictorially it is represented by the following grafting of trees, where the root of $\nu$ is grafted onto the $i$th leaf of $\mu$:
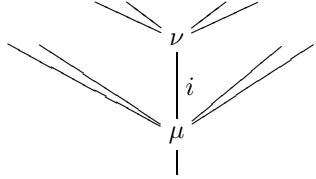
$$\xymatrix@R=6pt@C=6pt{*{} & *{} & *{} & *{} & *{} &  *{} & *{} & *{} & *{} & *{} &   *{} & *{} & *{} & \\
&&&&&&& \ar@{-}[ull] \ar@{-}[ul] \ar@{-}[ur] \ar@{-}[urr] \nu &    &&& &&\\
&&&&& && \quad i &&& &&& \\
&&&&&   &&\ar@{-}[uulllll] \ar@{-}[uullll] \ar@{-}[uu]\mu\ar@{-}[uurrr]  \ar@{-}[uurrrr] &&&  &&&\\
&&&&&   &&*{}\ar@{-}[u]&&&   &&&}$$

A ns operad can be defined as a graded vector space equipped with an identity operation and partial compositions which satisfy the unital axioms and the following two axioms (which replace associativity). For any $\lambda \in \PP_{l}, \mu \in \PP_{m}$ and $\nu \in \PP_{n}$ :

\medskip

Axiom I (sequential composition): 
$$\begin{array}{rcll}
(\lambda \circ_i \mu)\circ_{i-1+j}\nu  &=& \lambda \circ_i (\mu\circ_{j}\nu ), & 1\leq i\leq l, 1\leq j\leq m. 
\end{array}$$

\medskip
It corresponds to the two possibilities of composing in the following diagram:
$$
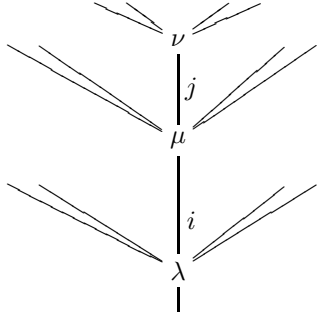
\xymatrix@R=6pt@C=6pt{
*{}&*{}&*{}&  *{}&*{}&*{}&*{}&*{}&   *{}&*{}&*{}&\\
&&&&&\ar@{-}[ull] \ar@{-}[ul] \ar@{-}[ur] \ar@{-}[urr]\nu&    &&& &&\\
&&& &&\quad j&&& &&& \\
&&&   &&\ar@{-}[uulllll] \ar@{-}[uullll] \ar@{-}[uu]\mu\ar@{-}[uurrr]  \ar@{-}[uurrrr] &&&  &&&\\
&&&   &&*{}\ar@{-}[u]&&&   &&&\\
&&& &&\quad i&&& &&& \\
&&&   &&\ar@{-}[uulllll] \ar@{-}[uullll] \ar@{-}[uu]\lambda\ar@{-}[uurrr]  \ar@{-}[uurrrr] &&&  &&&\\
&&&   &&*{}\ar@{-}[u]&&&   &&&
}$$

Axiom II (parallel composition):
$$\begin{array}{rcll}
(\lambda \circ_i \mu)\circ_{ k+m-1}\nu  &=& (\lambda \circ_k \nu)\circ_{i}\mu , &1\leq  i < k\leq l.
\end{array}$$

\medskip

It corresponds to the two possibilities of composing in the following diagram:
$$
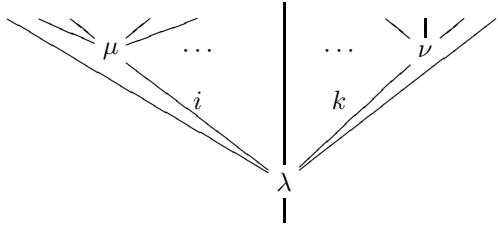
\xymatrix@R=6pt@C=6pt{
&&&&& &&&&& && \\
&&&&& &&&&& && \\
*{}&*{}&*{}&*{}&*{}&  *{}&*{}&*{}&*{}&*{}&   *{}&*{}&*{} \\
&&&\ar@{-}[ull] \ar@{-}[ul] \ar@{-}[ur] \ar@{-}[urr]\mu&&\cdots     && &\cdots &&
\ar@{-}[ul] \ar@{-}[u] \ar@{-}[ur] \nu&&\\
&&&&&i &&&k&& && \\
&&&&& &&& && && \\
&&&&&   &&\ar@{-}[uuuulllllll]\ar@{-}[uuullll]  \ar@{-}[uuuuu] \ar@{-}[uuurrr]\ar@{-}[uuuurrrrr]\lambda &&&  &&\\
&&&&&   &&*{}\ar@{-}[u]&&&   &&
}$$

\subsection{Symmetric operads}  A \emph{symmetric operad}, or operad for short, is a family of right $\Sy_{n}$-modules $\{\PP(n)\}_{n\geq 0}$ equipped with transformation of functors $\cc : \PP\circ \PP \to \PP$ and $\iota : \I \to \PP$ which are associative and unital. Here $\PP$ stands for the so-called \emph{Schur functor} defined as
$$\PP(V):= \bigoplus_{n}\PP(n)\t_{\Sy_{n}}V^{\t n}.$$

Any operad and any ns operad gives rise to a notion of algebras over this operad. There is a definition of operad using the partial compositions. See loc.\ cit.\ for details.  

\section{Associative algebras with a derivation}\label{S:preLiedend}

\subsection{Definition}\label{derivation} Let $A$ be nonunital associative algebra over $\KK$. A \emph{derivation} of $A$ is a linear map $D^{A}: A \to A$ which satisfies the Leibniz relation
$$D^{A}(ab) = D^{A}(a) b + a D^{A}(b)$$
for any $a,b\in A$. Let $(A', D^{A'})$ be another associative algebra with derivation. A morphism $f:(A,D^{A}) \to (A',D^{A'})$ is a linear map $f: A \to A'$ which is a morphism of associative algebras and which commutes with the derivations:
$$f\circ D^{A} = D^{A'} \circ f.$$

The operad governing the category of associative algebras with derivation admits the following presentation. There are two generating operations, one of arity 1, that we denote by $D$, and one of arity 2 that we denote by $\mu$. The relations are:
\begin{displaymath}
\left\{
\begin{array}{rcl}
\mu \circ (\mu, \id) & = & \mu \circ (\id, \mu),\\
D \circ \mu & = & \mu \circ (D, \id) + \mu \circ (\id, D).
\end{array} \right.
\end{displaymath}
They account for the associativity of the product and for the Leibniz relation. Since, in the relations, the variables stay in the same order, the category of associative algebras with derivation can be encoded by a nonsymmetric operad that we denote by $AsDer$. So it is determined by a certain family of vector spaces $AsDer_{n}, n\geq 1$, and composition maps
$$\cc_{i_{1}, \ldots , i_{k}} : AsDer_{k}\t AsDer_{i_{1}}\t \cdots \t AsDer_{i_{k}}\to AsDer_{n}$$
where $n = i_{1}+\cdots +i_{k}$.

\begin{thm}\label{AsDerexplicit} As a vector space $AsDer_{n}$ is isomorphic to the space of polynomials in $n$ variables:
$$AsDer_{n} = \KK[x_{1}, \ldots , x_{n}].$$
The composition map $\cc=\cc_{i_{1}, \ldots , i_{k}}$ is given by
\begin{displaymath}
\begin{array}{l}
\cc(P; Q_{1}, \ldots , Q_{k})(x_{1}, \ldots , x_{n}) =\\
P(x_{1}+\cdots +x_{i_{1}}, x_{i_{1}+1}+\cdots +x_{i_{1}+i_{2}}, x_{i_{1}+i_{2}+1}+ \cdots  , \ldots )Q_{1}(x_{1}, \ldots, x_{i_{1}}) Q_{2}(x_{i_{1}+1},\ldots ) \cdots \ .
\end{array}
\end{displaymath}
Under this identification the operations $\id, D, \mu$ correspond to $1_{1}, x_{1}\in \KK[x_{1}]$ and to $1_{2}\in \KK[x_{1}, x_{2}]$ respectively. More generally the operation 
$$(a_{1}, \ldots , a_{n})\mapsto D^{j_{1}}(a_{1})D^{j_{2}}(a_{2})\cdots D^{j_{n}}(a_{n})$$
 corresponds to the monomial $x_{1}^{j_{1}}x_{2}^{j_{2}}\cdots x_{n}^{j_{n}}$.
\end{thm}

Graphically the operation $x_{1}^{j_{1}}x_{2}^{j_{2}}\cdots x_{n}^{j_{n}}$ is pictured as a planar decorated tree as follows:

$$
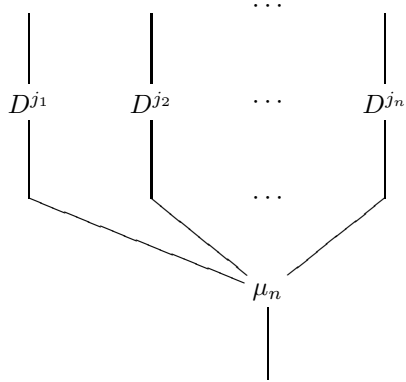
\xymatrix{
\ar@{-}[d] & \ar@{-}[d] & \cdots &  \ar@{-}[d] \\
D^{j_{1}}\ar@{-}[d]  & D^{j_{2}}\ar@{-}[d]  & \cdots & D^{j_{n}}\ar@{-}[d]  \\
*{}\ar@{-}[drr] & *{}\ar@{-}[dr] & \cdots  & *{} \ar@{-}[dl] \\
                    &                   &\mu_{n}\ar@{-}[d]    &   \\
                                &                   &*{}      &     }
$$
\begin{proo} Since $\mu$ is an associative operation the space of $n$-ary operations generated by $\mu$ is one-dimensional. Let us denote by $\mu_{n}$ the composition of $n-1$ copies of $\mu$. Using the Leibniz relation we see that any composite of copies of $\mu$ and $D$ can be written uniquely as composites of copies of $D$ first and then composites of copies of $\mu$. Hence any $n$-ary operation is  a linear combination of operations of the form 
$$\mu_{n}\circ (D^{j_{1}}, \ldots , D^{j_{n}})$$
for some sequence of nonnegative integers $\{j_{1}, \ldots, j_{n}\}$. We denote this operation by $x_{1}^{j_{1}}x_{2}^{j_{2}}\cdots x_{n}^{j_{n}}$. We obtain: $AsDer_{n} = \KK[x_{1}, \ldots , x_{n}]$.

In order to prove the formula for the composition map $\cc$, it is sufficient to prove that
\displaylines{
(P\circ_{i}Q)(x_{1}, \ldots , x_{n+m-1})= \hfill \cr
\qquad \qquad \qquad P(x_{1}, \ldots, x_{i-1}, x_{i}+\cdots x_{i+m-1},x_{i+m},  \ldots , x_{n+m-1})Q(x_{i}, \ldots , x_{i+m-1}),}
\noindent where $-\circ_{i}-$ is the $i$th partial composition, $P\in AsDer_{n},\  Q\in AsDer_{m}$.
It is sufficient to prove this formula when $P$ and $Q$ are monomials:

$$
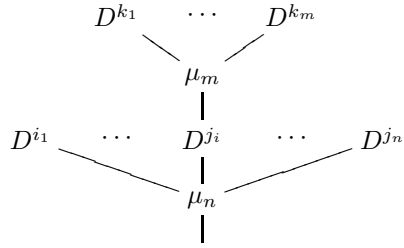
\xymatrix@R=10pt@C=10pt{
     & D^{k_{1}} & \cdots & D^{k_{m}} & \\
      &               &\mu_{m}\ar@{-}[ul]\ar@{-}[ur]\ar@{-}[d]&&\\
      D^{i_{1}} & \cdots & D^{j_{i}}  & \cdots & D^{j_{n}} \\
       & & \mu_{n}\ar@{-}[ull]\ar@{-}[urr]\ar@{-}[u]\ar@{-}[d]\\
       &&&& }$$

By direct inspection we see that it is sufficient to treat the case $Q= \mu_{m}$, and in fact the case $P= D^{\ell}, Q= \mu_{m}$. In other words we need to compute the element $D^{\ell}(a_{1}\cdots a_{m})$. By the Leibniz relation, for $\ell=1$, this is exactly the action of the operation $x_{1}+ \cdots +x_{m}$. Recursively we get that the operation $(a_{1}, \ldots , a_{m})\mapsto D^{\ell}(a_{1}\cdots a_{m})$ is 
$(x_{1}+ \cdots +x_{m})^{\ell}$ as expected. So we are done.
\end{proo}

\begin{lemma}\label{newdialg} The two binary operations $a \r b$ and $a \l b$ of $AsDer_{2}$ given by $x_{1}$ and $x_{2}$ respectively satisfy the following relations
\begin{displaymath}
\begin{array}{c}
(a \r b)\l c = a \r (b\l c), \\
(a \l b)\r c + (a \l b)\l c  = a \l (b\r c) +a \r (b\r c).
\end{array}
\end{displaymath}
\end{lemma}

\begin{proo} By a

straightforward operadic calculus we get
\begin{displaymath}
\begin{array}{cccc c}
x_{1}\circ_{1}x_{1}= (x_{1}+x_{2})x_{1} ,& x_{1}\circ_{1}x_{2}= (x_{1}+x_{2})x_{2} ,& x_{2}\circ_{1}x_{1}=  x_{3}x_{1} 
,& x_{2}\circ_{1}x_{2}= x_{3}x_{2} ,
\end{array}
\end{displaymath}
\begin{displaymath}
\begin{array}{cccc c}
x_{1}\circ_{2}x_{1}= x_{1}x_{2} ,& x_{1}\circ_{2}x_{2}= x_{1}x_{3} ,& x_{2}\circ_{2}x_{1}=  (x_{2}+x_{3})x_{2} 
,& x_{2}\circ_{2}x_{2}= (x_{2}+x_{3})x_{3}. \\
\end{array}
\end{displaymath}
It follows immediately that
\begin{displaymath}
\begin{array}{c}
x_{2}\circ_{1}x_{1}= x_{1}\circ_{2}x_{2}, \\
x_{1}\circ_{1}x_{2} +x_{2}\circ_{1}x_{2}= x_{1}\circ_{2}x_{1} + x_{2}\circ_{2}x_{1} .
\end{array}
\end{displaymath}
If we put $a\r b := x_{1}(a,b)$ and $a\l b := x_{2}(a,b)$, then we get the expected formulas.
\end{proo}

\subsection{Remarks} The relations of Lemma \ref{newdialg} give rise to a new type of  ns operad generated by two operations. Other examples include magmatic, dendriform \cite{JLLdig}, cubical \cite{LRhopf}, duplicial \cite{JLLgbto}, compatible-two-associative \cite{Dotsenko}.

Since $AsDer$ is a ns operad, it is completely determined by its free algebra over one generator. This free algebra has also been computed in \cite {GuoKeigher} by Guo and Keigher.

\subsection{An elementary example} Let $\KK=k[y]$ be the polynomial algebra in one variable over the field $k$. We consider the associative algebra with derivation $(A, D^{A})= (k[y][x], \frac {\partial}{\partial x})$. In $AsDer(k[y])^{\wedge}:= \Pi_{n}AsDer_{n}$ we consider the operation
$$exp(yD):= \sum_{n\geq 0}(\frac{y^n}{n!}D^n).$$
Then for any polynomial $p(x)\in k[x]$ the following formula holds in $k[y][x]= k[x,y]$:
$$exp(y\frac {\partial}{\partial x})(p(x))=p(x+y).$$
This is a key formula for studying vertex algebras, see for example \cite{Robinson}.

\section{First variation and dendriform algebras} Since the composition $\cc$ in the operad $AsDer$ needs only the sum, the product and the substitution of variables to be defined, it can be extended to many generalizations of the polynomial algebras: $C^{\infty}$-functions, rational functions, etc, provided that they are commutative. The rational function case is interesting since it permits us to treat integration, represented by the rational function $\frac{1_{1}}{x_{1}}$, which is the inverse of the operation $x_{1}$ coding for derivation. 

\subsection{The rational functions operad} We define a new ns operad $RatFct$ by
$$RatFct_{n}:= \KK( x_{1}, \ldots , x_{n})$$
(rational functions in the variables $x_{1}, \ldots , x_{n}$). The composition $\cc$ is given by the formula of Theorem \ref{AsDerexplicit}, or, equivalently by the partial composition formula:
\displaylines{
(P\circ_{i}Q)(x_{1}, \ldots , x_{n+m-1}):= \hfill \cr
\qquad \qquad \qquad P(x_{1}, \ldots, x_{i-1}, x_{i}+\cdots +x_{i+m-1},x_{i+m},  \ldots , x_{n+m-1})Q(x_{i}, \ldots , x_{i+m-1}).}

\begin{prop} The partial compositions $-\circ_{i}-$ as defined above make $RatFct$ into a ns operad.
\end{prop}
\begin{proo} There are two axioms to check, cf.\ \ref{nsop}. The first one follows from the fact that addition of variables is a formal group law, that is $F(x,y):= x+y$ satisfies $F(F(x,y), z) = F(x, F(y,z))$. The second axiom follows from the fact that the algebra of rational functions is commutative.
\end{proo}

We observe that the structure of associative algebra of $RatFct_{1}$ is precisely the algebra structure of the rational functions in one variable $\KK(x_{1})$.

\subsection{Integro-differential calculus} In the operad $RatFct$ the derivation operation $D$, represented by $x_{1}\in \KK(x_{1})=RatFct_{1}$, admits an inverse for composition, that is $1/x_{1}$, which is the integration operation $\int$ . So we can write the integro-differential calculus within the operad $RatFct$. Here is an example.

For integrable functions $f$ and $g$ on $\RR$ define 
$$(f\g g) (x) := f(x) \int_0^x g(t) dt \quad \mathrm{ and }\quad  (f\d g) (x) := \big(\int_0^x f(t) dt\big)\ g(x).$$
 Then, it is  shown in  \cite{JLLdig}  that these two operations satisfy the dendriform axioms (see below) as a consequence of integration by parts:
$$  \int_0^x g(t) dt  \int_0^x h(t) dt =  \int_0^x  g(t)\Big( \int_0^t h(s) ds\Big)dt +  \int_0^x \Big( \int_0^t g(s) ds\Big) h(t)dt\ .$$
We should be able to recover this property by computing in the operad $RatFct$ since $f\d g = (1/x_{1})(f,g)$ and $f\g g = (1/x_{2})(f,g)$, for $1/x_{1}, 1/x_{2}\in RatFct_{2}$. This is the object of Proposition \ref{dendtoratfct}.

More analogous formulas can be found in \cite{CHNT}.

\subsection{Dendriform algebras} Let us recall from \cite{JLLnote, JLLdig} that a \emph{dendriform algebra} is a vector space equipped with two binary operations $a\g b$ and $a\d b$ satisfying the relations 
\begin{displaymath}
\left \{ \begin{array}{rcl}
(a\g b) \g c & = & a\g (b\g c + b\d c),\\
(a\d b) \g c & = & a\d (b\g c),\\
(a\g b + a\d b) \d c & = & a\d (b\d c).
\end{array} \right.
\end{displaymath} 

It is proved in loc.\ cit.\ that the ns operad $Dend$, which encodes the dendriform algebras, is spanned in arity $n$ by the planar binary trees with $n+1$ leaves, a set that we denote by $PBT_{n+1}$. Let $Dend(\KK)= \bigoplus_{n} \KK[PBT_{n+1}]$ be the free dendriform algebra on one generator $\texttt{Y}$. The following relation holds in $Dend(\KK)$:
$$ s\vee t = s \d \texttt{Y} \g t\ ,$$
where $s\vee t$ stands for the grafting of the two binary trees $s$ and $t$.

\begin{prop}\label{dendtoratfct} There is an inclusion $\varphi$ of the ns operad $Dend$ to the ns operad $RatFct$ of rational functions induced by
\begin{displaymath}
\left \{ \begin{array}{rcl}
\d & \mapsto & \frac{1}{x_{1}}\ ,\\
\g & \mapsto & \frac{1}{x_{2}}\ . \\
\end{array} \right.
\end{displaymath}
It sends a planar binary tree $t$ to a rational function $\varphi(t)$ according to the following rules: $\varphi(|)=1$ and 
$$ \varphi(s\vee t) = \frac{\varphi(s)}{(x_{1}+\cdots + x_{p})}  \frac{\varphi(t)}{(x_{p+2}+\cdots + x_{p+q+1})}\in RatFct_{p+q+1} ,$$
when $s$, resp. \ $t$, is a tree with $p+1$ leaves, resp.\ $q+1$ leaves. 
\end{prop} 
\begin{proo} The fact that the two operations satisfy the dendriform axioms is a consequence of a 
straightforward operadic calculus:
\begin{displaymath}
\begin{array}{cccc c}
\frac{1}{x_{1}}\circ_{1}\frac{1}{x_{1}}= \frac{1}{x_{1}+x_{2}}\frac{1}{x_{1}} ,& \frac{1}{x_{1}}\circ_{1}\frac{1}{x_{2}}= \frac{1}{x_{1}+x_{2}}\frac{1}{x_{2}} ,& \frac{1}{x_{2}}\circ_{1}\frac{1}{x_{1}}= \frac{1}{x_{3}}\frac{1}{x_{1}}
,& \frac{1}{x_{2}}\circ_{1}\frac{1}{x_{2}}= \frac{1}{x_{3}}\frac{1}{x_{2}} ,
\end{array}
\end{displaymath}
\begin{displaymath}
\begin{array}{cccc c}
\frac{1}{x_{1}}\circ_{2}\frac{1}{x_{1}}=   \frac{1}{x_{1}}\frac{1}{x_{2}} ,& \frac{1}{x_{1}}\circ_{2}\frac{1}{x_{2}}=\frac{1}{x_{1}}\frac{1}{x_{3}} ,& \frac{1}{x_{2}}\circ_{2}\frac{1}{x_{1}}=  \frac{1}{x_{2}+x_{3}}\frac{1}{x_{2}}
,& \frac{1}{x_{2}}\circ_{2}\frac{1}{x_{2}}= \frac{1}{x_{2}+x_{3}}\frac{1}{x_{3}}. \\
\end{array}
\end{displaymath}

It follows immediately that
\begin{displaymath}
\begin{array}{c}
\frac{1}{x_{1}}\circ_{1}\frac{1}{x_{1}} +  \frac{1}{x_{1}}\circ_{1}\frac{1}{x_{2}}= \frac{1}{x_{1}+x_{2}}\frac{1}{x_{1}}+ \frac{1}{x_{1}+x_{2}}\frac{1}{x_{2}} =    \frac{1}{x_{1}}\frac{1}{x_{2}}=  \frac{1}{x_{1}}\circ_{2}\frac{1}{x_{1}} , \\
\\
 \frac{1}{x_{2}}\circ_{1}\frac{1}{x_{1}}= \frac{1}{x_{3}}\frac{1}{x_{1}}= \frac{1}{x_{1}}\circ_{2}\frac{1}{x_{2}} , \\
 \\
 \frac{1}{x_{2}}\circ_{1}\frac{1}{x_{2}} = \frac{1}{x_{3}}\frac{1}{x_{2}}= =  \frac{1}{x_{2}+x_{3}}\frac{1}{x_{2}} +  \frac{1}{x_{2}+x_{3}}\frac{1}{x_{3}}  =  \frac{1}{x_{2}}\circ_{2}\frac{1}{x_{1}} +  \frac{1}{x_{2}}\circ_{2}\frac{1}{x_{2}}  .
\end{array}
\end{displaymath}
If we put $a\d b := \frac{1}{x_{1}}(a,b)$ and $a\g b := \frac{1}{x_{2}}(a,b)$, then we get the expected formulas.

As for the second assertion, we first prove the formula for the grafting of trees. The proof is a straightforward operadic calculus, whose steps are the following. In $Dend(\KK)$ we have $s\vee t= s\d \texttt{Y} \g t= (s \d \texttt{Y}) \g t$. We first compute $ s\d \texttt{Y}= \cc(\d; s, \texttt{Y})= (\d\circ_{1}s)\circ_{p+1} \texttt{Y}$. Applying $\varphi$ we get
\begin{eqnarray*}
\varphi(s\d \texttt{Y})&=& \Big(\frac{1_{2}}{x_{1}}\circ_{1}s\Big)\circ_{p+1}\texttt{Y}\\
&=& \frac{1}{x_{1}+\cdots +x_{p}} \varphi(s)(x_{1}, \ldots , x_{p})\circ_{n}\texttt{Y}\\
&= & \frac{1}{x_{1}+\cdots +x_{p}} \varphi(s)(x_{1}, \ldots , x_{p}).
\end{eqnarray*}
A similar computation for the left product leads to the expected formula.

The rational functions $\varphi(t)$ for $t\in PBT_{n}$ are linearly independent (proof by induction). Hence $\varphi : Dend_{n} = \KK[PBT_{n+1}] \to RatFct_{n}$ is injective for all $n$.

\end{proo}

\subsection{Examples} Here is the image of the pb trees under $\varphi$ in low dimension:

$$
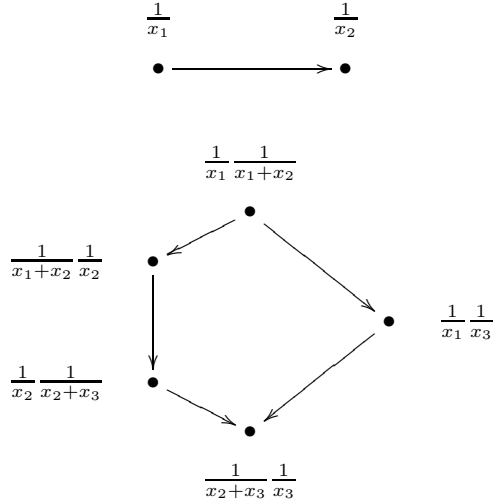
\xymatrix@R=4pt{
\frac{1}{x_{1}} & & \frac{1}{x_{2}}\\
\bullet\ar[rr] & &\bullet \\
}$$

$$\xymatrix@R=4pt@C=10pt{
 & & \frac{1}{x_{1}}\frac{1}{x_{1}+x_{2}}  &&& \\
 &&\bullet \ar[dl]\ar[ddrr]&&&\\
\frac{1}{x_{1}+x_{2}}  \frac{1}{x_{2}} & \bullet\ar[dd] & & && \\
& &&&\bullet\ar[ddll]& \frac{1}{x_{1}}\frac{1}{x_{3}}  \\
\frac{1}{x_{2}}\frac{1}{x_{2}+x_{3}} &\bullet \ar[dr] &&&&\\
&&\bullet&&&\\
&&\frac{1}{x_{2}+x_{3}}\frac{1}{x_{3}} &&&\\
}$$

\subsection{Associativity of $*$} Let us adopt the notation  $\Phi(u,v) = \frac{1}{u}+  \frac{1}{v}$. For any variables $x_{1}, x_{2}, x_{3}$ the formula 
$$\Phi(x_{1},x_{2})\Phi(x_{1}+x_{2}, x_{3})= \Phi(x_{2},x_{3})\Phi(x_{1},x_{2}+ x_{3})$$
is immediate to satisfy (cf. the proof of Proposition \ref{dendtoratfct}):
$$(\frac{1}{x_{1}}+\frac{1}{x_{2}})(\frac{1}{x_{1}+x_{2}}+\frac{1}{x_{3}})=\frac{x_{1}+ x_{2}+ x_{3}}{x_{1}x_{2}x_{3}}=
(\frac{1}{x_{2}}+\frac{1}{x_{3}})(\frac{1}{x_{1}}+\frac{1}{x_{2}+x_{3}}).$$
Viewed as an equality in the space $RatFct_{3}$ of ternary operations,  it simply says that the operation $* = \d + \g$ is associative. Indeed the left part of the equality is $\Phi\circ_{1}\Phi$ and the right part is $\Phi\circ_{2}\Phi$.

This formula is reminiscent of the cocycle condition of a 3-cochain $\phi$  in group cohomology:
$$\phi(x,y) + \phi(xy,z) = {}^x\phi(y,z)+\phi(x,yz),$$
and to the construction of the McLane invariant of a crossed module.

\subsection{Comparison with the operad of moulds} The work of J.\ Ecalle led F.\ Chapoton to introduce the \emph{operad of moulds} in \cite{Chapoton2, CHNT}. It is a nonsymmetric operad denoted by  $Mould$, which is determined by
$$Mould_{n}:= \KK( x_{1}, \ldots , x_{n}) ,$$
and by the partial compositions given by the formula

\displaylines{
(P\underline{\circ}_{i}Q)(x_{1}, \ldots , x_{n+m-1}):= \hfill \cr
( x_{i}+\cdots +x_{i+m-1})P(x_{1}, \ldots, x_{i-1}, x_{i}+\cdots +x_{i+m-1},x_{i+m},  \ldots , x_{n+m-1})Q(x_{i}, \ldots , x_{i+m-1}).} 

\begin{prop}\label{mould=ratfct} There is an isomorphism of ns operads $RatFract \cong Mould$.
\end{prop}

\begin{proo} By direct inspection we verify that the map $RatFract_{n}\to Mould_{n}$ given by
$$P(x_{1}, \ldots , x_{n})\mapsto (x_{1}+\cdots +x_{n})P(x_{1}, \ldots , x_{n})$$
is compatible with the operadic compositions. Since we are working with rational functions, the element $x_{1}+\cdots +x_{n}$ is invertible and so this map is an isomorphism.
\end{proo}

\subsection{Remarks} 

(a) Under the isomorphism of Proposition \ref{mould=ratfct}, Proposition \ref{dendtoratfct} can be found in  \cite{Chapoton2}. 

(b) Observe that the polynomials do not form a suboperad of $Mould$.

(c) If we think of the operation $1_{2}\in   \KK( x_{1},   x_{2})=RatFct_{2}$ as a third binary operation, then it is easy to check that the three operations $\d =\frac{1}{x_{1}}\ ,\g  =\frac{1}{x_{2}}\ , \cdot = 1_{2}$ satisfy the 7 axioms of a graded tridendriform algebra (cf.\ \cite{Chapoton1, LRtri}). We come back to this point in the next section.

\section{Second variation of $AsDer$}\label{secondvariation}

In probability theory there is a variation of the integration by parts called the \emph{Ito integral}. Its counterpart in the derivation framework, that could be called Ito derivation, is a linear map $D^{A}:A\to A$ which satisfies the following relation:
$$D^{A}(ab)= D^{A}(a) b + a D^{A}(b) + D^{A}(a) D^{A}(b).$$
One can treat both the derivation case, the Ito derivation case (and even more) by introducing a parameter $\lambda\in \KK$ as follows (we could also work with a \emph{formal} parameter $q$, that is take $\KK=k[q]$). 

\subsection{Parametrized derivation} By definition a \emph{$\lambda$-derivation} is a linear map $D^{A}:A\to A$, where $A$ is an associative algebra, such that 
$$D^{A}(ab)= D^{A}(a) b + a D^{A}(b) + \lambda D^{A}(a) D^{A}(b)$$
for any $a,b\in A$. For $\lambda=0$ we get the derivation, for  $\lambda=1$ we get the Ito derivation. By convention $\lambda=\infty$ stands for the case where $D^{A}$ is an algebra homomorphim: $D^{A}(ab)=D^{A}(a)D^{A}(b)$.

We denote by $\lambda\textrm{-}AsDer$ the operad of associative algebras equipped with a $\lambda$-derivation. 

We introduce the notation
$$\theta^{\lambda}(x_{1}, \ldots , x_{n}) := (x_{1}+\cdots +x_{n}) +\cdots +\lambda^{k-1} \theta_{k}(x_{1}, \ldots , x_{n})+\cdots +\lambda^{n-1}(x_{1}\ldots x_{n})$$
where $\theta_{k}(x_{1}, \ldots , x_{n})$ is the $k$th symmetric function of the variables $x_{1}, \ldots , x_{n}$. 

\begin{thm}\label{lambdaAsDerexplicit} As a vector space $\lambda\textrm{-}AsDer_{n}$ is isomorphic to the space of polynomials in $n$ variables:
$$\lambda\textrm{-}AsDer_{n} = \KK[x_{1}, \ldots , x_{n}].$$
The composition map $\cc$ is given by
\begin{displaymath}
\begin{array}{l}
\cc(P; Q_{1}, \ldots , Q_{k})(x_{1}, \ldots , x_{n}) =\\
{}\quad P(\theta ^{\lambda}(x_{1},\ldots ,x_{i_{1}}), \theta ^{\lambda}(x_{i_{1}+1},\ldots ,x_{i_{1}+i_{2}}) , \cdots )Q_{1}(x_{1}, \ldots x_{i_{1}}) Q_{2}(x_{i_{1}+i_{2}},\ldots ) \cdots \ .
\end{array}
\end{displaymath}
Under this identification the operations $\id, D, \mu$ correspond to $1_{1}, x_{1}\in \KK[x_{1}]$ and to $1_{2}\in \KK[x_{1}, x_{2}]$ respectively. More generally the operation 
$$(a_{1}, \ldots , a_{n})\mapsto D^{j_{1}}(a_{1})D^{j_{2}}(a_{2})\cdots D^{j_{n}}(a_{n})$$
 corresponds to the monomial $x_{1}^{j_{1}}x_{2}^{j_{2}}\cdots x_{n}^{j_{n}}$.
\end{thm}

\begin{proo} The proof is the same as in the case $
\lambda = 0$ performed in the first section. See also the proof of Proposition \ref{fgloperad}.
\end{proo}

\subsection{Remark} The associativity property of the composition in the operad $\lambda\textrm{-}AsDer$ implies that
$$\theta ^{\lambda}(x_{1},\ldots ,x_{i}, \theta ^{\lambda}(x_{i+1},\ldots ,x_{i+j}) , x_{i+j+1}, \ldots , x_{n}) = 
\theta ^{\lambda}(x_{1},\ldots ,x_{n}).$$
This formula can also be proved directly by observing that
$$1 + \lambda \theta ^{\lambda}(x_{1},\ldots ,x_{m})= \Pi_{i=1}^{m}(1+\lambda x_{i}).$$

\subsection{The parametrized operad $\lambda\textrm{-}RatFct$} As in section 2 we can put an operad structure on the rational functions by using the formulas of Theorem \ref{lambdaAsDerexplicit}. It gives a new operad, denoted  $\lambda\textrm{-}RatFct$, for which 
$$\lambda\textrm{-}RatFct_{n}= \KK( x_{1},\ldots , x_{n}),$$
and the partial composition is given by
$$\displaylines{
(P\circ_{i}Q)(x_{1}, \ldots , x_{n+m-1}):= \hfill \cr
\qquad \qquad P(x_{1}, \ldots, x_{i-1},\theta ^{\lambda}(x_{i},\cdots ,x_{i+m-1}),x_{i+m},  \ldots , x_{n+m-1})Q(x_{i}, \ldots , x_{i+m-1}).}
$$

\subsection{$\lambda$-TriDendriform algebras} In \cite{LRtri} we introduced the notion of tridendriform algebra, which is an algebra with 3 binary operations satisfying 7 relations (one for each of the cells of a triangle). The graded version was studied by Chapoton in \cite{Chapoton1}. There exists a parametrized version which handles both versions (and more) as follows. By definition a $\lambda$-tridendriform algebra has 3 binary operations denoted by $a\g b, a\d b, a\cdot b$ and 7 relations (one for each cell of the triangle):
\begin{displaymath}
\left\{
\begin{array}{lcr}
(x \g y) \g z &=& x \g (y * z)\ , \\
(x \d y) \g z &=& x \d (y \g z)\ , \\
(x * y) \d z &=& x \d (y \d z)\ , 
\end{array}\right.
\end{displaymath}
\begin{displaymath}
\left\{
\begin{array}{lcr}
(x \d y) \cdot z &=& x \d (y \cdot z)\ , \\ 
(x \g y) \cdot z &=& x \cdot (y \d z)\ , \\ (x \cdot y) \g z &=& x \cdot (y \g z)\ , \\ 
\end{array}\right.
\end{displaymath}
\begin{displaymath}
\left\{
\begin{array}{lcr}
(x \cdot y) \cdot z &=& x \cdot (y \cdot z)\ ,
\end{array}\right.
\end{displaymath}
where $x*y := x\g y + x\d y  + \lambda\  x\cdot y$. 

For $\lambda = 1$ we get Loday-Ronco's tridendriform algebra \cite{LRtri}, for  $\lambda = 0$ we get Chapoton's graded tridendriform algebra \cite{Chapoton1}.

\begin{prop} In the operad $\lambda\textrm{-}RatFct$ the binary operations
$$\d :=\frac{1_2}{x_{1}},\quad \g := \frac{1_2}{x_{2}},\quad \textrm{ and }\quad  \cdot := 1_{2}$$
satisfy the $\lambda$-tridendriform axioms.
\end{prop}

\begin{proo} The proof is an easy operadic computation analogous to the one performed in the proof of Proposition  \ref{dendtoratfct}. We leave it to the reader.
\end{proo}

\subsection{Formal group laws} Let $F(x,y)$ be a formal group law. This is a power series in the commutative variables $x$ and $y$ with coefficients in $\KK$ which satisfies the relations
$$F(x,0)= x = F(0, x) \quad , \quad F(F(x,y),z)=F(x,F(y,z)).$$
For instance, for any $\lambda \in \KK, \theta^{\lambda}(x,y):= x+y+\lambda xy$ is a polynomial formal group law. For $\lambda = 0$ it is the \emph{additive formal group law}, for $\lambda = 1$ it is the \emph{multiplicative formal group law}. Here is another example which is not polynomial. Let $\epsilon$ ad $\dd$ be elements in the commutative ring $\KK$ and let $Q(x):= 1- 2\dd x^2 + \epsilon x^4$. We suppose that 2 is invertible in $\KK$ so that $\sqrt{Q(x)}$ exists as a formal power series in $x$. The theory of elliptic curves gives rise to the \emph{Jacobi formal group law} defined by
$$ F(x,y) = \frac{ x\sqrt{Q(y)} + y\sqrt{Q(x)}}{1-\epsilon x^2 y^2}.$$

Let $F$ be any formal group law. By induction we define 
$$F(x_{1},\ldots, x_{n}) :=F(F(x_{1},\ldots, x_{n-1}), x_{n}) .$$

\begin{prop}\label{fgloperad} For any formal group law $F$ there is a well-defined nonsymmetric operad $RatFct^F$ given by:
$$RatFct^F_{n}:= \KK(x_{1}, \ldots ,x_{n}),$$
$$\displaylines{
(P\circ_{i}Q)(x_{1}, \ldots , x_{n+m-1}):= \hfill \cr
\qquad \qquad P(x_{1}, \ldots, x_{i-1},F(x_{i},\cdots ,x_{i+m-1}),x_{i+m},  \ldots , x_{n+m-1})Q(x_{i}, \ldots , x_{i+m-1}).}
$$\end{prop}
\begin{proo} Axiom I is an immediate consequence of the equality
$$F(x_{1},\ldots ,x_{i}, F(x_{i+1},\ldots ,x_{i+j}) , x_{i+j+1}, \ldots , x_{n}) = 
F(x_{1},\ldots ,x_{n}),$$
which follows from the associativity of $F$. 
Axiom II is a consequence of the commutativity of the algebra of rational functions.
\end{proo}

\section{$\PP$-algebras with derivation} Let $\PP$ be an algebraic operad, whose space of $n$-ary operations is denoted by $\PP(n)$. Recall that $\PP(n)$ is a right $\Sy_{n}$-module. A derivation on a $\PP$-algebra $A$ is a linear map $D^{A}: A\to A$ such that, for any operation $\mu\in \PP(n)$ the following formula holds:
$$D^{A}(\mu(a_{1}, \ldots, a_{n})) = \sum_{i=1}^{n} \mu(a_{1}, \ldots, D^{A}(a_{i}), \ldots, a_{n}).$$
 A similar computation as in the previous sections shows that the operad governing $\PP$-algebras with derivation, denoted by $\PP Der$, is such that 
 $$\PP Der(n)= \KK[x_{1}, \ldots, x_{n}]\t \PP(n)$$
 the action of the symmetric group $\Sy_{n}$ being the diagonal action (recall that $\Sy_{n}$ is acting on $ \KK[x_{1}, \ldots, x_{n}]$ by permuting the variables). The composition $\cc$ is obtained by combining the composition in $\PP$ and the formula in Theorem \ref{AsDerexplicit}.
 
 For instance $ComDer(n)= \KK[x_{1}, \ldots, x_{n}]$ is an $\Sy_{n}$-module and the composition map $\gamma$ of the operad $ComDer$ is given by the same formula as for the operad $AsDer$. It means that these formulas are compatible with the symmetric group action. Another way of phrasing this result is the following. Consider the forgetful functor which associates to a symmetric operad $\PP$ the nonsymmetric operad $\widetilde{\PP}$ such that $\widetilde{\PP}_{n}=\PP(n)$. Then we have $\widetilde{ComDer}=AsDer$.
 
 If we start with a ns operad $\PP$, then $\PP Der$ is a ns operad, where  $\PP Der_{n}= \KK[x_{1}, \ldots, x_{n}]\t \PP_{n}$. In terms of standard constructions in the operad framework, it is the Hadamard product of $\PP$ with $AsDer$.

\section{The ``pre-shuffle algebra'' of associative algebras with a derivation}

If we replace the space of polynomials (or formal power series) by the space of noncommutative polynomials (or noncommutative formal power series), then we do not get an operad anymore. However there is a variation of the notion of operads which permits us to provide a similar treatment in this noncommutative framework, it is the notion of ``pre-shuffle algebra'' due to M.\ Ronco \cite{Ronco07}
.

\subsection{Pre-shuffle algebra \cite{Ronco07}} A \emph{pre-shuffle algebra} 
is a family of vector spaces $\PP_{n}$ equipped with composition maps
$$\bullet_{i}: \PP_{m}\t \PP_{n} \to \PP_{m+n-1},\quad  1\leq i\leq m,$$
defined for $n\geq 1, m\geq 1$, which satisfy the following relations:
 $$(\lambda\bullet _{i}\mu)\bullet _{i-1+j }\nu= \lambda\bullet _{i}(\mu\bullet _{j}\nu)\ ,$$
for  $1\leq i\leq l, 1\leq j\leq m$. In other words, the difference with ns operads is that we keep only axiom I and we disregard axiom II (see \ref{nsop}), therefore any algebraic ns operad is a pre-shuffle algebra.

We denote by $PerAsDer$ the pre-shuffle algebra which is generated by a unary operation $D$ and a binary operation $\mu$, which satisfy the following relations:
\begin{displaymath}
\left\{
\begin{array}{rcl}
\mu \bullet_{1} \mu & = & \mu \bullet_{2} \mu\ ,\\
D \bullet_{1}\mu & = & \mu \bullet_{1}D + \mu \bullet_{2} D\ ,\\
(\alpha\bullet _{i}D)\bullet _{j}\mu &=& (\alpha\bullet _{j}\mu)\bullet _{i}D\ ,\\
(\alpha\bullet _{i}\mu)\bullet _{j+1}D &=& (\alpha\bullet _{j}D)\bullet _{i}\mu\ ,
\end{array} \right.
\end{displaymath}
for any operation $\alpha$ and $i<j$.

Observe that the first relation is the associativity of $\mu$, the second relation is saying that $D$ is a derivation, the third and fourth relations say that the operations $D$ and $\mu$ commute for parallel composition.

\begin{thm}\label{PerAsDerexplicit} As a vector space $PerAsDer_{n}$ is isomorphic to the space of noncommutative polynomials in $n$ variables:
$$PerAsDer_{n} = \KK\langle x_{1}, \ldots , x_{n}\rangle .$$
The composition map $\bullet_{i}$ is given by
\displaylines{
(P\bullet_{i}Q)(x_{1}, \ldots , x_{n+m-1})= \hfill \cr
\qquad \qquad \qquad P(x_{1}, \ldots, x_{i-1}, x_{i}+\cdots +x_{i+m-1},x_{i+m},  \ldots , x_{n+m-1})Q(x_{i}, \ldots , x_{i+m-1}).}

Under this identification the operations $\id, D, \mu$ correspond to $1_{1}, x_{1}\in \KK\langle x_{1} \rangle$ and to $1_{2}\in \KK\langle x_{1}, x_{2}\rangle$ respectively. More generally the operation 
$$\big(\cdot ((\mu\bullet_{j_{k}}D)\bullet_{j_{k-1}}D)\cdots \bullet_{j_{1}}D\big)$$
 corresponds to the noncommutative monomial $x_{j_{k}}x_{j_{k-1}}\cdots x_{j_{1}}$.
\end{thm}
Graphically the operation $x_{j_{k}}x_{j_{k-1}}\cdots x_{j_{1}}$ is pictured as a  planar decorated tree with levels as follows (example: $x_{1}x_{2}x_{n}x_{2}$) :
$$
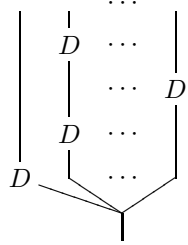
\xymatrix@R=4pt@C=4pt{
\ar@{-}[dddd] & \ar@{-}[d] & \cdots &  \ar@{-}[dd] \\
 & D\ar@{-}[dd]  & \cdots & \\
 &   & \cdots & D \ar@{-}[dd] \\
 & D\ar@{-}[d]  &  \cdots & \\
D\ar@{-}[drr]  &*{} \ar@{-}[dr]  & \cdots &*{} \ar@{-}[dl]  \\
                    &                   &*{} \ar@{-}[d]    &   \\
                                &                   &*{}      &     }
$$

\begin{proo} By \cite{Ronco07} we know that the free preshuffle algebra on a certain set of generating operations is spanned by some leveled planar trees whose vertices are labeled by the operations. Because of the relations entwining $D$ and $\mu$ we can move the operations $D$ up in a leveled tree composition of operations. Because of the associativity of $\mu$ the trees involving only $\mu$ give rise to corollas (as in the $AsDer$ case). Hence the operad $PerAsDer$ is spanned by the trees of the form indicated above. Note that the levels indicating the order in which the copies of the operation $D$ are performed is necessary in the preshuffle algebra framework.
\end{proo}

\subsection{Remark} One can also check that $PerAsDer$ is a shuffle algebra in the sense of \cite{Ronco07}.

\subsection{Formal group laws in noncommutative variables}  Let $F(x,y)$ be a formal group law in noncommutative variables. This is a series in $\KK((x,y))$ which satisfies the relation
$$F(F(x,y),z)=F(x,F(y,z)).$$
For instance, the  Baker-Campbell-Hausdorff series is the formal power series defined as
$$BCH(x,y):= \log(\exp(x)\exp(y)).$$
We recall that the first terms are 
$$BCH(x,y) = x+y+\frac{1}{2}[x,y] + \frac{1}{12}([[x,y],y]+[x,[x,y]]) + \cdots .$$

 By induction we define 
$$F(x_{1},\ldots, x_{n}) =F(F(x_{1},\ldots, x_{n-1}), x_{n}) .$$

\begin{prop} For any formal group law $F$ in noncommutative variables  there is a well-defined preshuffle algebra $PerRatFct^F$ given by:
$$PerRatFct^F_{n}:= \KK((x_{1}, \ldots ,x_{n})),$$
$$\displaylines{
(P\circ_{i}Q)(x_{1}, \ldots , x_{n+m-1}):= \hfill \cr
\qquad \qquad P(x_{1}, \ldots, x_{i-1},F(x_{i},\cdots ,x_{i+m-1}),x_{i+m},  \ldots , x_{n+m-1})Q(x_{i}, \ldots , x_{i+m-1}).}$$
If $F$ is polynomial, then the restriction to $\KK\langle x_{1}, \ldots , x_{n}\rangle $ is still a preshuffle algebra.

\end{prop}
\begin{proo} Axiom I is an immediate consequence of the equality
$$F(x_{1},\ldots ,x_{i}, F(x_{i+1},\ldots ,x_{i+j}) , x_{i+j+1}, \ldots , x_{n}) = 
F(x_{1},\ldots ,x_{n}).$$
The last assertion is immediate.
\end{proo}

\begin{cor} Let $F$ be the additive formal group law $F(x,y)=x+y$. Then the associated polynomial preshuffle algebra is $PerAsDer$.
\end{cor}
 
\section{Appendix: Homotopy associative algebras with derivation} We know that homotopy associative algebras are $\Ai$-algebras as defined by Jim Stasheff in  \cite{Stasheff}. Our purpose is to describe the notion of homotopy associative algebras with a derivation, that is to unravel the operad $AsDer_{\infty}$. The solution is given by the Koszul duality theory of quadratic operads, see for instance \cite{LV} where the Ginzburg-Kapranov theory is extended to operads generated by binary \emph{and} unary operations. A quadratic operad $\PP$ admits a Koszul dual cooperad $\PP\ac$. The cobar construction over $\PP\ac$ is the operad of $\PP$-algebras up to homotopy (i.e.\ the minimal model $\PP_{\infty}:= \Omega\, \PP\ac$ of the operad $\PP$) whenever the Koszul complex of the operad $\PP$ is acyclic. In this appendix we compute the cooperad $AsDer\ac$ and its linear dual $AsDer^{!}$, we prove that the Koszul complex $(AsDer\ac \circ AsDer, \dd)$ is acyclic and we unravel the co!
 bar construction $AsDer_{\infty}:= \Omega\ AsDer\ac$. So we get a precise description of the notion of \emph{associative algebra with derivation up to homotopy}.

If the parameter $\lambda\in \KK$ is different from $0$, then the operad $\lambda$-$AsDer$ is not a quadratic operad since the term $D(a)D(b)$ needs three generating operations to be defined. So one needs new techniques to extend Koszul duality to this case, see \cite{MerkulovVallette}.  

\subsection{The $AsDer^{!}$-algebras} The relations defining the operad $AsDer$ are quadratic since each monomial involves only the composition of two operations. Hence $AsDer$ is suitable for applying the Koszul duality theory. We use the notations and results of \cite{LV}.

\begin{prop} The Koszul dual operad of $AsDer$ is the operad $AsDer^{!}$ generated by the unary operation $d$ and the binary operation $\mu$ which satisfy the following relations:
$$d\circ d = 0,\quad   d\circ \mu= \mu\circ (d, \id)=  \mu\circ (\id, d), \quad \mu \circ (\mu, \id) = \mu\circ (\id, \mu) .$$
\end{prop}

In other words an $AsDer^!$-algebra is an associative algebra $A$ equipped with a linear map $d$ such that $d^2=0$ and $d(ab)=d(a)b=ad(b)$.

\begin{proo} The operad $AsDer$ is generated by the graded vector space
$$E=(0, \KK\, D, \KK\, \mu, 0, \ldots ).$$
 The weight $2$ subspace of the free operad $\TTT(E)$, denoted $\TTT(E)^{(2)}$,  is spanned by the operations which are composite of two of the generating operations. It is of dimension $6$ with basis
$$ \mu\circ_{1} \mu,\quad \mu\circ_{2} \mu,\quad  \mu\circ_{1} D,\quad  \mu\circ_{2} D,\quad  D\circ_{1} \mu,\quad  D\circ_{1} D.$$
The subspace of relations $R$ is of dimension $2$ spanned by
$D\circ_{1}\mu - \mu\circ_{1} D -  \mu\circ_{2} D $ and $ \mu\circ_{1} \mu -  \mu\circ_{2} \mu$. The Koszul dual cooperad $AsDer\ac$ is cogenerated by $sD$ and $s\mu$ ($s$ is the shift of degree), with $s^2 R$ as corelations: 
$$AsDer\ac= \Id \oplus sE \oplus s^2R\oplus \cdots .$$
By definition the Koszul dual operad $AsDer^!$ of $AsDer$ is, essentially, the linear dual of $AsDer\ac$. Let us denote by $d$ the linear dual of $sD$ (put in degree 0) and $\mu$ the linear dual of $s\mu$ (put in degree 0). Then the space 
$ \TTT(\KK d \oplus \KK \mu)^{(2)}$ is also of dimension 6 and the quotient
 $AsDer^!_{2}=  R^{\vee} = \TTT(\KK d \oplus \KK \mu)/ R^{\perp}$ is two dimensional. So $R^{\perp}$ is the 4-dimensional space spanned by the elements
$$d\circ_{1}d,\quad  d\circ_{1}\mu - \mu\circ_{1}d,\quad  d\circ_{1}\mu - \mu\circ_{2}d,\quad  \mu\circ_{1}\mu- \mu\circ_{2}\mu.$$

\end{proo}

\begin{prop}\label{AsDerduale} The space of $n$-ary operations of the operad $AsDer^{!}$ is 2-dimensional:
$$AsDer^{!}_{n}= \KK\,\mu_{n}\oplus  \KK\,d\mu_{n}, n\geq 2, \quad and \quad  AsDer^{!}_{1}= \KK\,\id\oplus  \KK\, d .$$
The partial composition $\circ_{i}$ is given by
\begin{displaymath}
\begin{array}{rcl}
\mu_{m}\circ_i \mu_{n}&=& \mu_{m+n-1},\\
 d\mu_{m}\circ_i \mu_{n}&=& d\mu_{m+n-1},\\
  \mu_{m}\circ_i d\mu_{n}&=& d\mu_{m+n-1},\\
   d\mu_{m}\circ_i d\mu_{n}&=& 0,
\end{array}
\end{displaymath}

where by convention $\mu_{1}=\id$ (so $d\circ_{1} \mu_{n} = d\mu_{n}$).
\end{prop}
\begin{proo} The generating binary operation $\mu$ generates the operation $\mu_{n}$ in arity $n$. The relations entwining $\mu$ and $d$ imply that the only other possibility to create an operation in arity $n$ is to compose with a copy of $d$. The formula for the partial composition is obtained by direct inspection.
\end{proo}

\begin{prop}\label{AsDeracyclic} The operad $AsDer$ is a Koszul operad.
\end{prop}
\begin{proo} Most of the methods for proving Koszul duality would work in this simple case. We choose to write down explicitly the ``rewriting system method'', see \cite{LV}.

The operad $AsDer$ is presented by the generators $D$ and $\mu$ and the rewriting relations
\begin{displaymath}
\left\{
\begin{array}{rcl}
\mu\circ_{1}\mu &\mapsto& \mu\circ_{2}\mu ,\\
D\circ_{1}\mu &\mapsto& \mu\circ_{1}D + \mu\circ_{2}D.
\end{array} \right.
\end{displaymath}

The critical monomials are $\mu\circ_{1}(\mu\circ_{1}\mu)$ and $D\circ_{1}(\mu\circ_{1}\mu)$. The first critical monomial is known to be confluent (Koszulity of the operad $As$), but let us recall the proof. One one hand, one has
$$((xy)z)t \mapsto (x(yz))t \mapsto x((yz)t) \mapsto x(y(zt)).$$
On the other hand, one has 
$$((xy)z)t \mapsto (xy)(zt) \mapsto x(y(zt)).$$
Since one ends up with the same element, we have shown that the first critical monomial is confluent.

Let us show confluence for the second critical monomial.  One one hand, one has
$\begin{array}{l}D((xy)z) \mapsto D(x(yz)) \mapsto(Dx)(yz) \mapsto x(D(yz)) \mapsto (Dx)(yz)+ x((Dy)z+y(Dz)) \mapsto  \end{array}$

$\begin{array}{r} (Dx)(yz)+ x((Dy)z)+x(y(Dz)). \end{array}$

On the other hand, one has 

\noindent $\begin{array}{l}D((xy)z) \mapsto (D(xy))z) + (xy)(Dz) \mapsto ((Dx)y)z +(x(Dy))z + (xy)(Dz) \mapsto \end{array}$

$\begin{array}{r}(Dx)(yz)+ x((Dy)z)+x(y(Dz)).\end{array}$

Since one ends up with the same element, we are done.

Since the critical monomials are confluent the ns operad $AsDer$ is Koszul.
\end{proo}

\subsection{Homotopy associative algebra with derivation} Since the operad $AsDer$ is a Koszul operad, its minimal model is given by the operad $\Omega AsDer\ac$. The aim of this section is to describe this operad explicitly.

By definition a \emph{homotopy associative algebra with derivation} is an algebra over the differential graded ns operad  $AsDer_{\infty}$ constructed as follows. In arity $n$ the space $(AsDer_{\infty})_{n}$ is spanned by the planar trees with $n$ leaves whose nodes with $k$ inputs are labelled by $m_{k}$ or $Dm_{k}$ when $k\geq 2$ and by $D$ when $k=1$. Observe that the symbol 
$Dm_{k}$ is to be taken as a whole and not as $m_{k}$ followed by $D$. The homological degree of $m_{k}$ and of $Dm_{k}$ is $k-1$ for $k\geq 2$, the homological degree of $D$ is 0.

Example:
$$
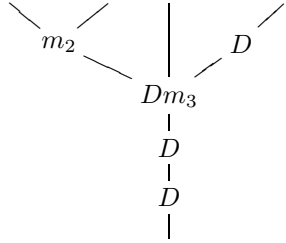
\xymatrix@R=6pt@C=6pt{
*{}\ar@{-}[dr] & &*{}\ar@{-}[dl]  & *{}\ar@{-}[dd] & & *{}\ar@{-}[dl]  \\
& m_{2}\ar@{-}[drr]  & & &D\ar@{-}[dl]  &  \\
& & &Dm_{3} \ar@{-}[d]  & &  \\
& & & D  \ar@{-}[d]  & &  \\
& & & D  \ar@{-}[d]  & &  \\
& & &*{}& &  
}$$

The differential map $\partial : (AsDer_{\infty})_{n}\to (AsDer_{\infty})_{n}$ of homological degree $-1$ is induced by
$$
\begin{array}{rcl}
\partial (D) & = & 0,\\
\partial(m_{n})&=& - \sum_{\substack{ n=p+q+r \\ k=p+1+r\\ k>1, q>1}} (-1)^{p+qr} m_{k}\circ (\id^{\t p}\t m_{q} \t \id^{\t r}), \\
\partial(Dm_{n})&=& - \sum_{\substack{ n=p+q+r \\ k=p+1+r\\ k\geq 1, q> 1}} (-1)^{p+qr} Dm_{k}\circ (\id^{\t p}\t m_{q} \t \id^{\t r})\\
&& - \sum_{\substack{ n=p+q+r \\ k=p+1+r\\ k> 1, q\geq 1}} (-1)^{p+qr} m_{k}\circ (\id^{\t p}\t Dm_{q} \t \id^{\t r})
\end{array}
$$
for $n\geq 2$.

For instance, in low arity we get:
$$
\begin{array}{rcl}
\partial (m_{2}) & = & 0,\\
\partial(m_{3})&=&m_{2}\circ (\id, m_{2}) - m_{2}\circ (m_{2}, \id), \\
\partial(Dm_{2})&=&- D\circ m_{2} + m_{2}\circ (\id, D) + m_{2}\circ (D, \id).
\end{array}
$$
\begin{prop} The operad $AsDer_{\infty}$ is isomorphic to the minimal model $\Omega AsDer\ac$ of the operad $AsDer$. Hence a homotopy  associative algebra with derivation is an algebra over the operad  $AsDer_{\infty}$.
\end{prop}

\begin{proo} By definition the cobar construction over the cooperad $AsDer_{\infty}\ac$ is the free operad $\TTT(s \overline {AsDer\ac})$. Since, in arity $n\geq 2$, the space $(AsDer\ac)_{n}$ is spanned by two elements and the space $ \overline {AsDer\ac}_{1}$ by one element (the notation overline means that we get rid of the identity operation), the free operad $\TTT(s \overline {AsDer\ac})$ is spanned by the planar rooted trees with labelled nodes as in the description of $AsDer_{\infty}$. In the free operad the operadic composition is given by grafting.

The operad $\Omega AsDer\ac$ is a differential graded operad, so we need to describe the differential map. It is sufficient to describe it on the operadic generators, that is on the corollas labelled by either $m_{k}$ or $Dm_{k}$. This boundary map is deduced from the cooperad structure of $AsDer\ac$, that is, from the operad structure of $AsDer^{!}$. From the formulas in Proposition \ref{AsDerduale} we deduce the formulas given in the construction of the operad   $AsDer_{\infty}$.
\end{proo}

\subsection{Transfer theorem} Let us recall that the interest of the notion of ``algebra up to homotopy'' lies, in part, in the following transfer theorem. Let $(A,\dd )$ be a differential graded associative algebra with derivation $D^{A}$ of degree $0$, that is $\dd (ab)= \dd (a)b + (-1)^{|a|}a \dd (b)$ and $\dd (D^{A}(a)) = D^{A}(\dd (a))$, and let $(V,\dd )$ be a retract by deformation of the chain complex $(A,\dd )$ (e.g.\ $(H(A), 0)$ when $\KK$ is a field). Then $(V,\dd )$ can be equipped with a $AsDer_{\infty}$-algebra structure transferred from the $AsDer$-algebra structure of $(A,\dd )$. In particular there exist analogues of the Massey products on $H(A)$. Since $AsDer$ is Koszul, all these results are particular examples of \cite{LV}, Chapter 9, which extend the results obtained by Tornike Kadeishvili \cite{Kadeishvili} on differential graded associative algebras.


\begin{thebibliography}{99}

\bibitem{Chapoton1} \emph{F.\ Chapoton}, Alg\`ebres de Hopf des permuto\`edres, associa\`edres et hypercubes,
Advances in Math. (150), no 2 (2000), 264--275.

\bibitem{Chapoton2} \emph{F.\ Chapoton}, The anticyclic operad of moulds,
International Math. Research Notices (2007) no 20.

\bibitem{CHNT} \emph{F.\ Chapoton, F.\ Hivert, J.-C.\ Novelli, J.-Y.\ Thibon}, An operational calculus for the Mould operad, Int. Math. Res. Not. IMRN 2008, no. 9, Art. ID rnn018, 22 pp.

\bibitem{Dotsenko} \emph{V.\ Dotsenko}, Compatible associative products and trees, J.\ Algebra \& Number Theory, 3 (2009), no. 5, 567-586.

\bibitem{DotsenkoKhoroshkin} \emph{V.\ Dotsenko, A.\ Khoroshkin}, Gr\"obner basis for operads, 	 Duke Math.\ Journal, to appear.

\bibitem{Kadeishvili} \emph{T.\ Kadeishvili,} The algebraic structure in the homology of an $A(\infty)$-algebra. Soobshch. Akad. Nauk Gruzin. SSR, 108:249--252, 1982.

\bibitem{GuoKeigher} \emph{Li Guo; W.\ Keigher},
On differential Rota-Baxter algebras.
J.\ Pure Appl.\ Algebra 212 (2008), no. 3, 522--540. 
 
\bibitem{JLLnote} \emph{J.-L.\ Loday}, Alg\`ebres ayant deux op\' erations associatives (dig\`ebres). C.\ R.\ Acad.\ Sci.\ Paris S\'er. I Math. 321 (1995), no. 2, 141--146.

\bibitem{JLLdig} \emph{J.-L.\ Loday},  Dialgebras, in ``Dialgebras and related operads", Springer Lecture Notes in Math. 1763 (2001), 7-66.

\bibitem{JLLgbto} \emph{J.-L.\ Loday}, Generalized bialgebras and triples of operads, Ast\' erisque 320 (2008), x+116 pp.

\bibitem{LRhopf}  \emph{J.-L.\ Loday,  M.\ Ronco}, Hopf algebra of the planar binary trees. Adv. in Maths 139 (1998), 293--309.

\bibitem{LRtri}  \emph{J.-L.\ Loday,  M.\ Ronco}, Trialgebras and families of polytopes,  in "Homotopy Theory: Relations with Algebraic Geometry, Group Cohomology, and Algebraic K-theory" Contemporary Mathematics 346 (2004), 369--398.

 
 \bibitem{LV}  \emph{J.-L.\ Loday,  B.\ Vallette}, Algebraic operads, in preparation.
  
\bibitem{MSS}  \emph{M.\ Markl,  J.\ Stasheff, S.\ Shnider}, Operads in algebra, topology and physics. Mathematical Surveys and Monographs, 96. American Mathematical Society, Providence, RI, 2002. x+349 pp.

\bibitem{MerkulovVallette}   \emph{ S.\ Merkulov, B.\ Vallette}, Deformation theory of representations of prop(erad)s. I. J. Reine Angew. Math. 634 (2009), 51--106.

\bibitem{Robinson} \emph{T.\ Robinson},
New perspectives on exponentiated derivations, the formal 
Taylor theorem, and Fa\`a di Bruno's formula, arXiv math:0903.3991.

\bibitem{Ronco07} \emph{ M.\ Ronco},  Shuffle bialgebras, Ann.\ Inst.\ Fourier (2010), to appear.



\bibitem{Stasheff} \emph{J.D. \ Stasheff}, Homotopy associativity of $H$-spaces. I, II. Trans. Amer. Math. Soc. 108 (1963), 275-292; ibid. 293--312. 

\end{thebibliography}
\end{document}